\documentclass{amsart}
\usepackage{graphics}

\newtheorem{theorem}{Theorem}[section]
\newtheorem{proposition}[theorem]{Proposition}
\newtheorem{lemma}[theorem]{Lemma}
\newtheorem{claim}[theorem]{Claim}

\theoremstyle{definition}
\newtheorem{definition}[theorem]{Definition}

\theoremstyle{remark}
\newtheorem{remark}[theorem]{Remark}

\numberwithin{equation}{section}

\newcommand{\PL}{\mathbb{P}^1}
\newcommand{\AL}{\mathbb{A}^1}
\newcommand{\CT}{\mathbb{C}^\times}
\newcommand{\G}{\gamma}
\newcommand{\D}{\delta}
\newcommand{\Aut}{\hbox{Aut}}
\newcommand{\spc}{\hspace{.3cm}}

\begin{document}

\title{Mirror symmetry and $\CT$}
\author{Nobuyoshi Takahashi}
\address{Department of Mathematics, 
Hiroshima University, Higashi-Hiroshima 739-8526, Japan}
\email{takahasi@math.sci.hiroshima-u.ac.jp}


\begin{abstract}
We show that counting functions of covers of $\CT$ 
are equal to sums of integrals associated to certain `Feynman' graphs. 
This is an analogue of the mirror symmetry 
for elliptic curves(\cite{Dijkgraaf}).
\end{abstract}

\maketitle

\section{Introduction and main result}
According to Mirror Symmetry, 
it is believed that 
various counting functions on Calabi-Yau manifolds 
are related to period integrals on their mirror families. 
One case where the relation is well understood is 
that of 1-dimensional compact Calabi-Yau manifolds, 
that is, elliptic curves. 
In \cite{Dijkgraaf}, it is shown 
that the counting function of simply ramified finite coverings 
of an elliptic curve 
is equal to the partition functions 
in fermionic and bosonic field theories on the mirror family. 

We can also consider non-compact but algebraic analogues 
of Calabi-Yau manifolds. 
The simplest one is $\CT$: 
we regard it as an `open Calabi-Yau manifold' 
since it has a nowhere vanishing holomorphic 1-form $dz/z$ 
that has only logarithmic poles at $0$ and $\infty$. 
We consider the following enumerative problem 
and show that the resulting functions 
are equal to bosonic partition functions on $\CT$. 

\begin{definition}\label{n}
Let $b$ be a nonnegative integer 
and $k, l, d_1, \dots, d_k, e_1, \dots, e_l$ positive integers. 
We consider $\CT$ as $\PL\setminus\{0, \infty\}$ and 
take distinct points $P_1, \dots, P_b\in\CT$ 
and sufficiently small circles $\G$ and $\G'$ 
around $0$ and $\infty$ respectively. 
Let $\D$ (resp. $\D'$) be the disjoint union of 
$d_i$-ple (resp. $e_i$-ple) covers of $\G$ (resp. $\G'$), 
and let $p:\D\rightarrow \CT$ 
(resp. $p':\D'\rightarrow \CT$) be the natural map. 

Then we define 
$n_{b;d_1, \dots, d_k; e_1, \dots, e_l}$ to be the number 
of isomorphism classes of quadruples $(C, \pi, i, i')$ 
where:
\begin{itemize}
\item $C$ is a possibly non-connected smooth curve, 
\item $\pi:C\rightarrow\CT$ is a finite covering 
that is simply ramified over $P_1, \dots, P_b$ 
(i.e., exactly one branch maps with degree 2) and 
unramified elsewhere, and 
\item $i:\D\rightarrow C$ and $i':\D'\rightarrow C$ satisfy 
$\pi\circ i=p$ and $\pi\circ i'=p'$. 
\end{itemize}

Here, an isomorphism 
of $(C_1, \pi_1, i_1, i'_1)$ and $(C_2, \pi_2, i_2, i'_2)$ is 
an isomorphism $f:C_1\to C_2$ such that 
$\pi_2\circ f=\pi_1$, $f\circ i_1=i_2$ and $f\circ i'_1=i'_2$. 
\end{definition}

\begin{definition}
Let $F_{b,k,l}(z_1, \dots, z_k; w_1, \dots, w_l)$ be 
the generating function 
\[
  \sum n_{b;d_1, \dots, d_k; e_1, \dots, e_l}
      z_1^{d_1}.\cdots.z_k^{d_k}.w_1^{-e_1}.\cdots.w_l^{-e_l}. 
\]
\end{definition}

\begin{remark}
(1) The number $n_{b;d_1, \dots, d_k; e_1, \dots, e_l}$ is nonzero 
only if $\sum d_i = \sum e_i$. 
Thus, $F_{b,k,l}$ is homogeneous of degree 0. 

(2) If $k+l+b$ is odd, we have $F_{b,k,l}=0$. 
\end{remark}

\begin{figure}
\caption{\label{graphs}}
\vspace{.5cm}
\parbox{\linewidth}{
\parbox[t]{.5\linewidth}
{$b=4, k=l=1$.

\vspace{.2cm}
\includegraphics{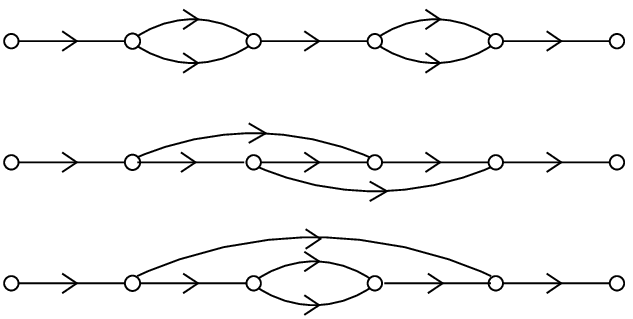}}
\parbox[t]{.5\linewidth}
{$b=3, k=2, l=1$. 

\vspace{.2cm}
\includegraphics{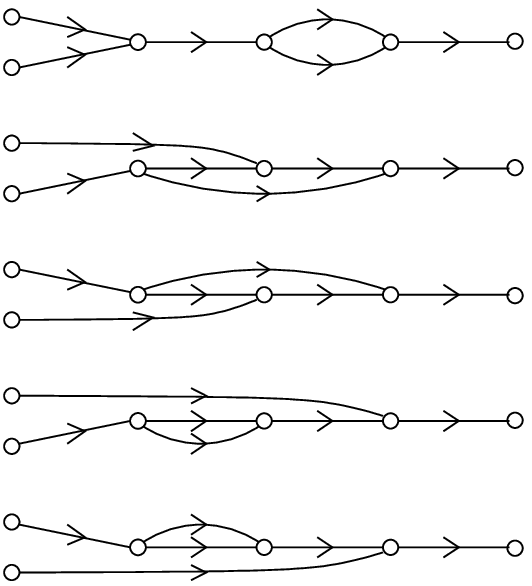}}
}
\end{figure}

For bosonic interpretation, 
we consider the following `Feynman' graphs. 
Figure \ref{graphs} shows examples. 

\begin{definition}\label{defgraph}
Let $b$ be a nonnegative integer and $k, l$ positive integers. 
We consider sets 
$V_1=\{z_1, \dots, z_k\}$, $V_2=\{x_1, \dots, x_b\}$, 
$V_3=\{w_1, \dots, w_l\}$ and $V=V_1\coprod V_2\coprod V_3$. 

We define a partial ordering on $V$ as follows: 
for $v, v'\in V$, $v<v'$ 
if and only if $v\in V_i$, $v'\in V_{i'}$ with $i<i'$ 
or $v=x_j, v'=x_{j'}$ and $j<j'$. 

Let $E$ be a set and $v_{in}, v_{fin}$ maps $E\to V$. 
Then $V$ and $E$ defines an oriented graph. 
Now we define $G_{b, k, l}$ as the set of graphs satisfying 
\begin{itemize}
\item for any $e\in E$, $v_{in}(e)<v_{fin}(e)$, 
\item $\#\{e\in E|v_{in}(e)=z_i\}=1$, 
\item  
$\{\#\{e\in E|v_{in}(e)=x_i\}, \#\{e\in E|v_{fin}(e)=x_i\}\}=\{1,2\}$, and 
\item $\#\{e\in E|v_{fin}(e)=w_i\}=1$. 
\end{itemize}

We define an isomorphism of graphs 
$(E, v_{in}, v_{fin})$ and $(E', v'_{in}, v'_{fin})$ 
to be a bijection of $E$ and $E'$ 
that commutes with $v_{in}, v_{fin}$ and $v'_{in}, v'_{fin}$. 
\end{definition}

The following is our `Bosons' formula. 

\begin{theorem}\label{boson}
In the notations of Definition \ref{defgraph}, 
we regard elements of $V$ as variables. 
Then, if $\max\{|z_i|\}< r_1 <\dots< r_b < \min\{|w_i|\}$, 
we have 
\[
F_{b,k,l}(z_1, \dots, z_k; w_1, \dots, w_l)=
\sum_{[\Gamma]\in G_{b,k,l}} \frac{I_\Gamma}{\#\Aut(\Gamma)}, 
\]
where, for $\Gamma=(E, v_{in}, v_{fin})$, 
$I_\Gamma$ is defined as 
\[
  \prod_{i=1}^b\left(\oint_{|x_i|=r_i}\frac{dx_i}{2\pi\sqrt{-1}x_i}\right)
    \prod_{e\in E}\frac{v_{in}(e)v_{fin}(e)}{(v_{in}(e)-v_{fin}(e))^2}. 
\]
\end{theorem}

\begin{remark}
Consider a graph as in Definition \ref{defgraph} with 
the condition 
$\{\#\{e\in E|v_{in}(e)=x_i\}, \#\{e\in E|v_{fin}(e)=x_i\}\}=\{1,2\}$
replaced by 
$\{\#\{e\in E|v_{in}(e)=x_i\}, \#\{e\in E|v_{fin}(e)=x_i\}\}=\{0,3\}$. 
For such a graph, the integral in Theorem \ref{boson} is 0 
since all the poles of the integrand with respect to $x_i$ 
are either all contained in $|x_i| < r_i$ or in $|x_i| > r_i$. 

Thus, $F_{b, k, l}$ is also equal to 
the sum of the integrals as in Theorem \ref{boson} 
over the set of isomorphism classes of graphs $\Gamma$ 
with `incoming', `intermediate' and `outgoing' vertices 
$z_1, \dots, z_k; x_1, \dots, x_b; w_1, \dots, w_l$ 
--- ordering of intermediate vertices is essential ---
satisfying the following conditions: 

\begin{itemize}
\item the graph has no loops, that is, 
the two end points of any edge are not the same vertex, 
\item each vertex $x_i$ has three edges connected, 
\item each $z_i$ is connected 
to exactly one of $x_1,\dots,x_b, w_1,\dots,w_l$ 
and each $w_i$ is connected 
to exactly one of $x_1,\dots,x_b, z_1,\dots,z_k$. 
\end{itemize}
\end{remark}

\begin{remark}
(1) 
As for `Fermions' formula in \cite{Dijkgraaf}, 
it is straightforward to modify the arguments in \S5 of \cite{Dijkgraaf} 
to prove the following formula:
\[
\sum_{b=0}^{\infty} F_{b,1,1}(q;1)\frac{\lambda^b}{b!} =  
    \left( \sum _{p\in \mathbb{Z}_{\geq 0}+1/2}
         q^{p}e^{\lambda p^{2}/2}\right) 
    \left( \sum _{p\in \mathbb {Z}_{\geq 0}+1/2}
         q^{p}e^{-\lambda p^{2}/2}\right).
\]
\cite{O} shows that the generating function for 
all $n_{b;d_1, \dots, d_k; e_1, \dots, e_l}$'s
satisfies the Toda lattice hierarchy of Ueno and Takasaki.

(2) The related problem 
of counting the numbers of covers 
of $\PL$ with arbitrary ramification over $\infty$, 
called Hurwitz numbers, 
is also being actively studied
(see, for example, \cite{SSV}, \cite{Vakil}, \cite{GJ1}, \cite{GJ2}). 
This can be seen as an analogue of enumeration of curves in Fano manifolds, 
and as in the case of Fano manifolds, 
some recursive formulas are known. 

(3) In dimension 2, the author counted curves of degree up to 8 
in an affine cubic surface 
which are images of morphisms from an affine line(\cite{Takahashi}), 
and interpreted those numbers 
from the view point of mirror symmetry of log surfaces(\cite{Takahashi2}). 
\end{remark}

\section{Proof}

Now we prove Theorem \ref{boson}. 
We use notations in Definitions \ref{n} and \ref{defgraph}. 

Let $d_1, \dots, d_k$ and $e_1, \dots, e_l$ be positive integers 
such that $\sum d_i=\sum e_i=d$. 
We rewrite $n_{b;d_1, \dots, d_k; e_1, \dots, e_l}$ 
as in \S4 of \cite{Dijkgraaf}. 

\begin{definition}
We denote the group of permutations of $\{1, \dots, d\}$ by $S_d$. 

Let $\tilde{d}_i = \sum_{j=1}^i d_j$, 
and let $\sigma_{d_1, \dots, d_k}\in S_d$ be the permutation 
\[
(\tilde{d}_0+1\spc\dots\spc\tilde{d}_1)
(\tilde{d}_1+1\spc\dots\spc \tilde{d}_2)
.\cdots.(\tilde{d}_{k-1}+1\spc\dots\spc\tilde{d}_k). 
\]
We define $\tilde{e}_i$ and $\sigma_{e_1, \dots, e_l}$
in the same way. 

For $\sigma, \tau\in S_d$, we write 
$\sigma^\tau=\tau^{-1}\sigma\tau=\tau\circ\sigma\circ\tau^{-1}$. 
\end{definition}

\begin{lemma}
$n_{b;d_1, \dots, d_k; e_1, \dots, e_l}$ 
is the number of sequences $(g_1, \dots, g_b, \tau)$, 
where $g_i\in S_d$ are transpositions and $\tau\in S_d$, such that 
$g_b.\cdots.g_1.\sigma_{d_1,\dots,d_k}=(\sigma_{e_1,\dots,e_l})^\tau$. 
\end{lemma}

\begin{proof}
We may take $\G=\{|z|=r_0\}$, $\G'=\{|z|=r_b\}$ 
and $P_i=x_i e^{\epsilon\sqrt{-1}}$, 
where $0 < r_0 < x_1 < \dots < x_b < r_b$ and $0 < \epsilon \ll 1$. 
Also choose $r_i$ with $x_i < r_i < x_{i+1}$ for $i=1, \dots, b-1$. 
Then, let $\G_i$ be the cycle that starts at $r_0$, 
goes to $r_i$, then runs on $|z|=r_i$ counter-clockwise 
and goes back to $r_0$(see Figure \ref{cycles}).

\begin{figure}
\caption{\label{cycles}}
\vspace{.5cm}
\includegraphics{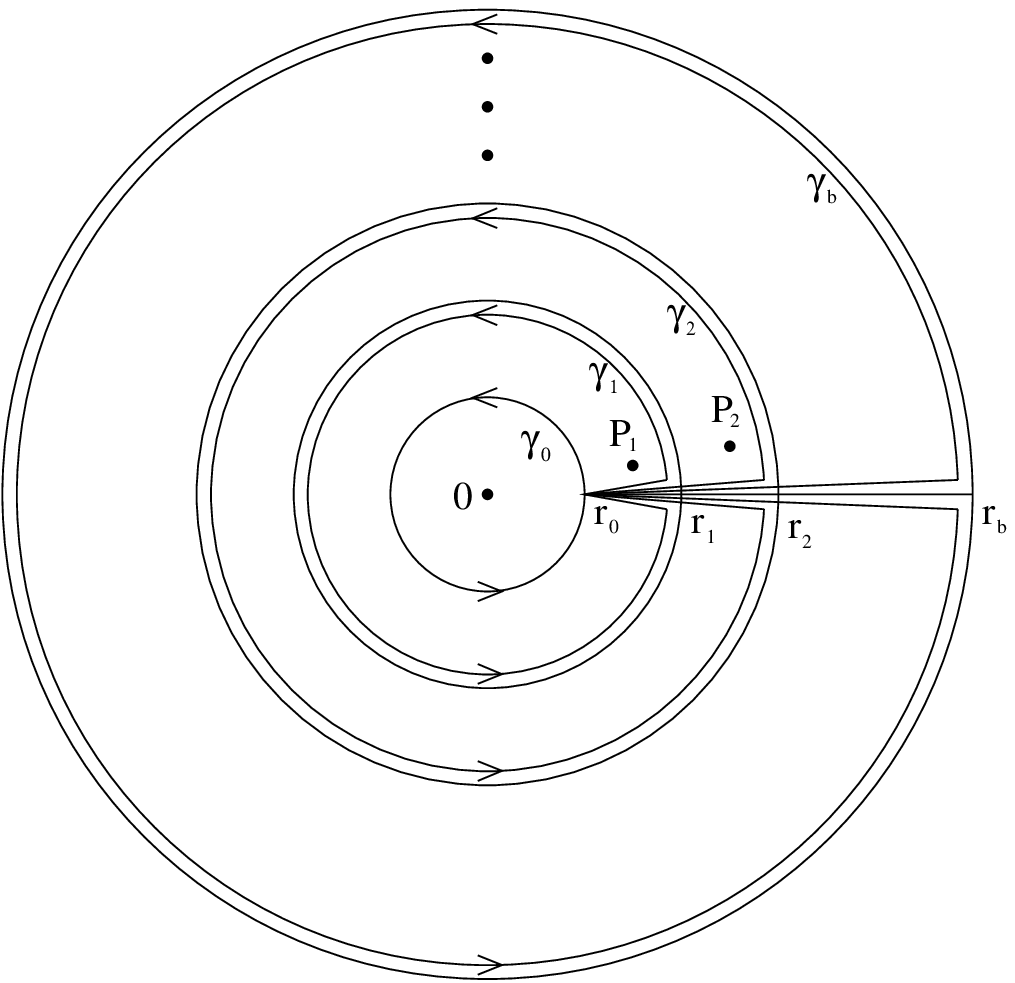}
\end{figure}

We choose a numbering of $p^{-1}(r_0)$, 
where $p:\D\to\CT$ is as in Definition \ref{n}, 
such that the monodromy of $p$ along $\G_0$ is 
$\sigma_{d_1, \dots, d_k}$. 
Similarly, for $p^{\prime-1}(r_b)$, 
we give a numbering so that the monodromy 
along the circle $|z|=r_b$, counter-clockwise, is 
$\sigma_{e_1, \dots, e_l}$. 

Let $(C, \pi, i, i')$ be given. 
Then, 
we can talk of monodromy $\tilde{g}_i\in S_d$ of $\pi$ along $\G_i$. 
Then $g_i:=\tilde{g}_i \tilde{g}_{i-1}^{-1}$ is a transposition, 
since it is the monodromy of the simple ramification over $P_i$. 
Also, let $\tau\in S_d$ be the permutation defined by 
the path from $r_b$ to $r_0$ 
and numbering of $p^{\prime-1}(r_b)$ and $p^{-1}(r_0)$. 
Then, we have 
$g_b.\cdots.g_1.\sigma_{d_1,\dots,d_k}
=\tilde{g}_b=(\sigma_{e_1,\dots,e_l})^\tau$. 

Conversely, such $(g_1, \dots, g_b, \tau)$ defines $(C, \pi, i, i')$. 
\end{proof}

We associate a graph to each $(g_1, \dots, g_b, \tau)$ 
as in the above lemma. 

\begin{definition}\label{assgraph}
Let $(g_1, \dots, g_b, \tau)$ be as in the above lemma. 
For $i=0, \dots, b$, let $E_i$ be the set of cyclic components of 
$g_i.\cdots.g_1\sigma_{d_1, \dots, d_k}$. 

Let $\sim$ be the equivalence relation on $\coprod E_i$ 
generated by the following relation: 
for $\sigma\in E_i$ and $\sigma'\in E_{i+1}$, 
$\sigma\sim\sigma'$ if $\sigma=\sigma'$ as elements of $S_d$. 

Then let $E = \coprod E_i/\sim$. 
For an element 
$e=\{\sigma_{i_0}, \sigma_{i_0+1}, \dots, \sigma_{i_1}\}$ of $E$, 
where $\sigma_i\in E_i$, 
we define $v_{in}(e)$ to be $x_{i_0}$ if $i_0>0$
and $v_{fin}(e)$ to be $x_{i_1+1}$ if $i_1<b$. 
If $i_0=0$, we take $v$ to satisfy 
$\sigma_0=(\tilde{d}_{v-1}+1\spc\dots\spc\tilde{d}_v)$
and let $v_{in}(e)=z_v$, 
and similarly if $i_1=b$, we take $v$ with
$\sigma_b=(\tilde{e}_{v-1}+1\spc\dots\spc\tilde{e}_v)^\tau$
and let $v_{fin}(e)=w_v$. 

We say $(g_1, \dots, g_b, \tau)$ belongs to 
(the isomorphism class of) $(E, v_{in}, v_{fin})$. 
\end{definition}

\begin{remark}
In the situation of the above lemma and definition, 
$E_i$ is naturally in one-to-one correspondence 
with the set of connected components of $\pi^{-1}(\{|z|=r_i\})$. 
The graph defined describes how these cycles 
meet and break. 
\end{remark}

\begin{lemma}
The isomorphism class of $(E, v_{in}, v_{fin})$ 
defined in Definition \ref{assgraph} 
belongs to $G_{b, k, l}$. 
\end{lemma}

\begin{proof}
We use notations in Definition \ref{assgraph}. 

If $\sigma_i\in E_i$ and $e$ is the equivalence class of $\sigma_i$, 
then $v_{in}(e)=x_i$ if and only if $g_i$ and $\sigma_i$ intersect 
and $v_{fin}(e)=x_{i+1}$ if and only if $\sigma_i$ and $g_{i+1}$ 
intersect. 

There are either one or two cyclic components of 
$g_{i-1}.\cdots.g_1\sigma_{d_1, \dots, d_k}$ that intersect $g_i$, 
since $g_i$ is a transposition. 
If there is only one, say $\sigma_{i-1}$, 
then two components of $g_i.\cdots.g_1\sigma_{d_1, \dots, d_k}$ 
intersect $g_i$, i.e. the two components of $g_i \sigma_{i-1}$. 
In this case, we have 
$\#\{e\in E|v_{in}(e)=x_i\}=2$ and 
$\#\{e\in E|v_{fin}(e)=x_i\}=1$. 
If there are two, $\sigma_{i-1}$ and $\sigma'_{i-1}$, 
then $g_i \sigma_{i-1}\sigma'_{i-1}$ is 
the sole component of $g_i.\cdots.g_1\sigma_{d_1, \dots, d_k}$ 
that intersects $g_i$, 
and we have $\#\{e\in E|v_{in}(e)=x_i\}=1$, 
$\#\{e\in E|v_{fin}(e)=x_i\}=2$. 

The other properties are easy to see. 
\end{proof}

\begin{definition}
For a graph $\Gamma=(E,v_{in},v_{fin})$, 
let $n_{\Gamma; d_1, \dots, d_k; e_1, \dots, e_l}$ 
be the number of sequences $(g_1, \dots, g_b, \tau)$, 
where $g_i$ are transpositions and $\tau\in S_d$, such that 
$g_b.\cdots.g_1.\sigma_{d_1, \dots, d_k}=(\sigma_{e_1, \dots, e_l})^\tau$ 
and that the associated graph is isomorphic to $\Gamma$. 

Also, let 
\[
F_\Gamma(z_1, \dots, z_k; w_1, \dots, w_l)=
  \sum
    n_{\Gamma; d_1, \dots, d_k; e_1, \dots, e_l}
      z_1^{d_1}.\cdots.z_k^{d_k}.w_1^{-e_1}.\cdots.w_l^{-e_l}. 
\]
\end{definition}

From the two lemmas above, 
it suffices to prove the following to prove Theorem \ref{boson}.  

\begin{proposition}\label{propboson}
Let $\Gamma$ be a graph whose isomorphism class belongs 
to $G_{b, k, l}$. 
Then, if $\max\{|z_i|\}<\min\{|w_i|\}$, we have 
\[
F_\Gamma(z_1, \dots, z_k; w_1, \dots, w_l)= \frac{I_\Gamma}{\#\Aut(\Gamma)}, 
\]
where, as in Theorem \ref{boson}, 
\[
 I_{\Gamma} := 
  \prod_{i=1}^b\left(\oint_{|x_i|=r_i}\frac{dx_i}{2\pi\sqrt{-1}x_i}\right)
    \prod_{e\in E}\frac{v_{in}(e)v_{fin}(e)}{(v_{in}(e)-v_{fin}(e))^2}, 
\]
with $\max\{|z_i|\}< r_1 <\dots< r_b < \min\{|w_i|\}$. 
\end{proposition}

We prove the proposition by induction on $b$. 

\begin{claim}
Proposition \ref{propboson} holds for $b=0$.  
\end{claim}
\begin{proof}
In this case, we have $k=l$ 
and $\Gamma$ connects $z_i$ to $w_{\sigma(i)}$ 
for some $\sigma\in S_k$. 
By symmetry, we may assume $\sigma=id$. 
Then, $\tau\in S_d$ satisfies 
$\sigma_{d_1, \dots, d_k}=(\sigma_{e_1, \dots, e_k})^\tau$
and the associated graph is isomorphic to $\Gamma$ 
if and only if $d_i=e_i$ and 
\[
(\tilde{d}_{i-1}+1\spc\dots\spc\tilde{d}_i) = 
(\tilde{e}_{i-1}+1\spc\dots\spc\tilde{e}_i)^\tau
\]
for $i=1, \dots, k$. 
When $d_i=e_i$, the number of such $\tau$ is $\prod d_i$. 
Now 
\[
\sum_{d_1, \dots, d_k=1}^\infty \prod d_i 
   \prod(z_i^{d_i}w_i^{-d_i})
 = \prod \frac{z_i w_i}{(z_i-w_i)^2} 
\]
proves the claim. 
\end{proof}

\begin{claim}
If Proposition \ref{propboson} holds for $b=n$, 
it holds for $b=n+1$. 
\end{claim}

\begin{proof}
Let $\Gamma=(E, v_{in}, v_{fin})$ be an element of $G_{n+1,k,l}$. 

Case (1): $\#\{e\in E|v_{in}(e)=x_{n+1}\}=1$. 

By symmetry, we may assume $v_{fin}(e)=w_l$. 
There are two edges $e', e''$ 
such that $v_{fin}(e')=v_{fin}(e'')=x_{n+1}$. 
Let $\Gamma'\in G_{n,k,l+1}$ be the graph obtained by 
removing $x_{n+1}$ 
and connecting $e'$ to $w_l$ and $e''$ to $w_{l+1}$. 

Write $\tilde{e}_{i}=\sum_{j=1}^i e_j$.
If $(g_1, \dots, g_{n+1}, \tau)$ belongs to $\Gamma$ 
and $g_{n+1}.\cdots.g_1\sigma_{d_1, \dots, d_k}
=(\sigma_{e_1, \dots, e_l})^\tau$, 
then there exist $\tau'\in S_d$ 
with $\tau'(i)=\tau(i)$ for $0\leq i\leq\tilde{e}_{l-1}$ 
and positive integers $e'_l, e'_{l+1}$ with $e'_l+e'_{l+1}=e_l$ 
such that $g_n.\cdots.g_1\sigma_{d_1, \dots, d_k}
=(\sigma_{e_1, \dots, e_{l-1}, e'_l, e'_{l+1}})^{\tau'}$ 
and that $(g_1, \dots, g_n, \tau')$ belongs to $\Gamma'$, 
as is easily seen from definition. 

Let $e_1, \dots, e_{l-1}, e'_l, e'_{l+1}$ be positive integers 
with $\sum_{i=1}^{l-1}e_i+e'_l+e'_{l+1}=d$. 
Then, let $X_{e'_l, e'_{l+1}}$ be the set of 
$(g_1, \dots, g_n, t_1, \dots, t_{\tilde{e}_{l-1}})$ 
such that there exists $\tau'\in S_d$ with the property that 
$\tau'(i)=t_i$($0\leq i\leq\tilde{e}_{l-1}$), 
$g_n.\cdots.g_1.\sigma_{d_1, \dots, d_k}=
(\sigma_{e_1, \dots, e_{l-1}, e'_l, e'_{l+1}})^{\tau'}$ 
and that the graph associated to $(g_1, \dots, g_n, \tau')$ 
is isomorphic to $\Gamma'$. 

If $v_{in}(e')\not=v_{in}(e'')$ or $e'_l\not=e'_{l+1}$, 
$\tau'$ maps $(\tilde{e}'_{l-1}+1\spc\dots\spc\tilde{e}'_l)$ 
to the cyclic component of $g_n.\cdots.g_1.\sigma_{d_1, \dots, d_k}$ 
connected to $v_{in}(e')$, 
and $(\tilde{e}'_l+1\spc\dots\spc\tilde{e}'_{l+1})$ 
to the one connected to $v_{in}(e'')$. 
If $e'_l\not=e'_{l+1}$, 
$\tau'$ maps $(\tilde{e}'_{l-1}+1\spc\dots\spc\tilde{e}'_l)$ 
and $(\tilde{e}'_l+1\spc\dots\spc\tilde{e}'_{l+1})$ 
to cyclic components of lengths $e'_l$ and $e'_{l+1}$, respectively. 
Thus there are $e'_l.e'_{l+1}$ choices for $\tau'$ in these cases 
and we have 
$n_{\Gamma'; d_1,\dots, d_k; e_1, \dots, e_{l-1}, e'_l, e'_{l+1}}
=\# X_{e'_l, e'_{l+1}}.e'_l.e'_{l+1}$. 
If $v_{in}(e')=v_{in}(e'')$ and $e'_l=e'_{l+1}$ hold, 
then it is $2\# X_{e'_l, e'_{l+1}}.e_l^{\prime 2}$.  

On the other hand, for each 
$(g_1, \dots, g_n, t_1, \dots, t_{\tilde{e}_{l-1}})\in X_{e'_l, e'_{l+1}}$, 
there are $e'_l e'_{l+1}$ transpositions $g_{n+1}$ 
for which there exists $\tau$ 
with $\tau(i)=t_i$ for $0\leq i\leq\tilde{e}_{l-1}$ 
such that $(g_1, \dots, g_{n+1}, \tau)$ belongs to $\Gamma$. 
The number of such $\tau$ is $e_l:=e'_l + e'_{l+1}$, and we have 
\[
n_{\Gamma; d_1, \dots, d_k; e_1, \dots, e_l}=
\sum \# X_{e'_l, e'_{l+1}}.e_l.e'_l.e'_{l+1}. 
\]
Here the sum is taken over \emph{different} $X_{e'_l, e'_{l+1}}$'s 
with $e'_l+e'_{l+1}=e_l$. 

Noting that $X_{e'_l, e'_{l+1}}$ and $X_{e'_{l+1}, e'_l}$ are the same 
for $e'_l\not= e'_{l+1}$ 
if and only if $v_{in}(e')=v_{in}(e'')$, we see 
\[
n_{\Gamma; d_1, \dots, d_k; e_1, \dots, e_l}=
\sum_{e'_l+e'_{l+1}=e_l}
e_l n_{\Gamma'; d_1, \dots, d_k; e_1, \dots, e_{l-1}, e'_l, e'_{l+1}}
\]
if $v_{in}(e')\not=v_{in}(e'')$ and 
\[
n_{\Gamma; d_1, \dots, d_k; e_1, \dots, e_l}=
\frac{1}{2}\sum_{e'_l+e'_{l+1}=e_l}
e_l n_{\Gamma'; d_1, \dots, d_k; e_1, \dots, e_{l-1}, e'_l, e'_{l+1}}
\]
if $v_{in}(e')=v_{in}(e'')$. 

Finally, since 
\[
I_\Gamma=
\oint\frac{dx}{2\pi\sqrt{-1}x}\frac{x w_l}{(x-w_l)^2}
I_{\Gamma'}(z_1, \dots, z_k; w_1, \dots, w_{l-1}, x, x)
\]
holds, 
\[
\oint\frac{dx}{2\pi\sqrt{-1}x}\frac{x w_l}{(x-w_l)^2}x^{-e'_l}x^{-e'_{l+1}}
  = e_l w_l^{-e_l}
\]
proves the assertion. 

Case (2): $\#\{e\in E|v_{in}(e)=x_{n+1}\}=2$. 

By symmetry, we may assume 
that $w_{l-1}$ and $w_l$ appear as $v_{fin}(e)$. 
Let $e'$ be the edge such that $v_{fin}(e')=x_{n+1}$. 
Let $\Gamma'\in G_{n,k,l-1}$ be the graph obtained by 
removing $x_{n+1}$ and connecting $e'$ to $w_{l-1}$. 

As in Case (1), 
If $(g_1, \dots, g_{n+1}, \tau)$ belongs to $\Gamma$ 
and $g_{n+1}.\cdots.g_1\sigma_{d_1, \dots, d_k}
=(\sigma_{e_1, \dots, e_l})^\tau$, 
then there exists $\tau'\in S_d$ 
such that $g_n.\cdots.g_1\sigma_{d_1, \dots, d_k}
=(\sigma_{e_1, \dots, e_{l-2}, e'_l})^{\tau'}$, 
where $e'_l=e_{l-1}+e_l$, 
and that $(g_1, \dots, g_n, \tau')$ belongs to $\Gamma'$. 

For positive integers $e_1, \dots, e_{l-2}, e'_{l-1}$ 
with $\sum_{i=1}^{l-2}e_i+e'_{l-1}=d$, 
let $X$ bet the set of 
$(g_1, \dots, g_n, t_1, \dots, t_{\tilde{e}_{l-2}})$ 
such that there exists $\tau'\in S_d$ with the property that 
$\tau'(i)=t_i$($0\leq i\leq\tilde{e}_{l-2}$), 
$g_n.\cdots.g_1.\sigma_{d_1, \dots, d_k}=
(\sigma_{e_1, \dots, e_{l-2}, e'_{l-1}})^{\tau'}$ 
and that the graph associated to $(g_1, \dots, g_n, \tau')$ 
is isomorphic to $\Gamma'$. 

Then there are $e'_{l-1}$ such $\tau'$, and 
$n_{\Gamma'; d_1,\dots, d_k; e_1, \dots, e_{l-2}, e'_{l-1}}$ is 
$\# X.e'_{l-1}$. 

On the other hand, 
for $e_{l-1}$ and $e_l$ with $e_{l-1}+e_l=e'_{l-1}$, 
the number of pairs $(g_{n+1}, \tau)$ 
with $\tau(i)=t_i$($0\leq i\leq\tilde{e}_{l-2}$) 
such that $(g_1, \dots, g_{n+1}, \tau)$ belongs to $\Gamma$ 
is $e'_{l-1}$ times $e_{l-1} e_l$: 
let $(a_1\spc\dots\spc a_{e'_{l-1}})$ be the cyclic permutation 
corresponding to $e'$. 
Choose $1\leq r\leq e'_{l-1}$ 
and let $s$ be the residue of $r+e_l$ modulo $e'_{l-1}$. 
We take $g_{n+1}$ to be $(a_r\spc a_s)$, 
and then $g_{n+1}.(a_1\spc\dots\spc a_{e'_{l-1}})$ has cyclic components 
of lengths $e_{l-1}$ and $e_l$, 
containing $a_r$ and $a_s$ respectively. 
We can choose $\tau$ so that 
the two components equal 
$(\tilde{e}_{l-2}+1\spc\dots\spc\tilde{e}_{l-1})^\tau$ and 
$(\tilde{e}_{l-1}+1\spc\dots\spc\tilde{e}_l)^\tau$, 
respectively. 

Thus we have 
\[
n_{\Gamma; d_1, \dots, d_k; e_1, \dots, e_l}=
\#X.e'_{l-1}.e_{l-1}.e_l = 
e_{l-1}e_l n_{\Gamma'; d_1, \dots, d_k; e_1, \dots, e'_{l-1}}. 
\]
Since 
\[
I_\Gamma=
\oint\frac{dx}{2\pi\sqrt{-1}x}
\frac{x w_{l-1}}{(x-w_{l-1})^2}\frac{x w_l}{(x-w_l)^2}
I_{\Gamma'}(z_1, \dots, z_k; w_1, \dots, w_{l-2}, x) 
\]
holds in this case, 
\[
\oint\frac{dx}{2\pi\sqrt{-1}x}
    \frac{x w_{l-1}}{(x-w_{l-1})^2}\frac{x w_l}{(x-w_l)^2}x^{-e'_{l-1}}
  = \sum_{e_{l-1}+e_l=e'_{l-1}} e_{l-1} e_l w_{l-1}^{-e_{l-1}} w_l^{-e_l}
\]
proves the assertion. 

\end{proof}

\end{document}